\documentclass[ a4paper,reqno]{amsart}

\usepackage{amsthm}
\usepackage{amsfonts}
\usepackage{amsmath}

\usepackage{mathptmx}      
 \newtheorem{theorem}{Theorem}[section]
 \newtheorem{corollary}[theorem]{Corollary}
 \newtheorem{lemma}[theorem]{Lemma}
 
 \theoremstyle{definition}
 
 \theoremstyle{remark}
 \newtheorem{remark}[theorem]{Remark}
 
 \numberwithin{equation}{section}

\def\ox{x}
\def\xn{x_{n+1}}
\def\xv{\overline{x}}
\def\uv{{{u}}}

\def\scal#1#2{\langle #1, #2\rangle }
\def\M{{\mathcal M}}

\def\R#1{{\mathbb R}^{#1}}

\def\Sph{\mathbf{S}^{n}}
\def\Sirc{\mathbf{S}^{n-1}}
\def\dsi{}
\def\Disk{\mathbf{B}^{n}}
\def\Ball{\mathbf{B}^{n+1}}
\def\diver{\mathrm{div}}
\begin{document}

\title{Disjoint minimal graphs}
\thanks{The paper was supported by grant of RFBR no.~03-01-00304.}


\author{Vladimir G. Tkachev }


\address{Volgograd State University\\\textit{Current address}: Royal Institute of Technology, Stockholm}

\email{tkatchev@kth.se}

\maketitle

\begin{abstract}
We prove that the number $s(n)$ of disjoint minimal graphs
supported on domains in $\R{n}$ is bounded by $e(n+1)^2$. In the
two-dimensional case we show that $s(2)\leq 3$.

\end{abstract}


%
%

\section{Introduction}

Let $w(x)$ be a solution to the minimal surface equation
\begin{equation}\label{mineq}
\diver \frac{\nabla w(x)}{\sqrt{1+|\nabla w(x)|^2}}=0
\end{equation}
defined in an open subset  $G\subset \R{n}$. The graph $\mathcal{G}=(G,w)$:
\begin{equation}\label{xn}
\xn=w(\ox), \qquad \ox\in G
\end{equation}
is called a minimal graph \textit{supported} on $G$ if $w(\ox)$
changes no sign on $G$ and
\begin{equation}\label{e:supp}
\quad w(\ox)=0,\quad \ox \in\partial G.
\end{equation}
In this case the domain $G$ is said to be \textit{admissible}; a
finite collection of admissible disjoint domains $\{G_j\}_{j=1}^{s}$
will also called  admissible. It follows from the maximum principle
for solutions of (\ref{mineq}), (\ref{e:supp}) that an admissible
domain is necessarily unbounded.

An admissible domain $G$ is called \textit{trivial} if the
complement $\R{2}\setminus \overline{G}$ contains no unbounded
components.  It follows immediately from the maximum principle that if an admissible collection
$\{G_j\}_{j=1}^s$ contains a trivial domain then $s=1$.

The following problem was recently posed in \cite{MR}: \textit{How
many admissible domains $G_j$ one can arrange in $\R{n}$
without overlapping}? P.~Li and J.~Wang \cite{Li0} proved that in
an arbitrary  dimension $n$ there are only finitely many minimal
graphs supported on disjoint open subsets. Let us denote by $s(n)$
the maximal cardinality $s$ of an admissible collection in
dimension $n$. In \cite{Li0} the following uniform estimate was
proven
\begin{equation}\label{double}
s(n)\leq 2^{n+1}(n+1).
\end{equation}
On the other hand, a conjecture of W.~Meeks \cite{MR} states that
$s(2)=2$. In \cite{Spruck}, J.~Spruck showed that this property
holds if minimal graphs have a
sublinear growth.

In this paper we obtain two effective estimates on $s(n)$ without
any additional requirements. More precisely, we will prove

\begin{theorem}\label{th:main3}
Let $n\geq 3$. Then the following (polynomial in $n$) estimate
holds
$$
s(n)\leq e(n+1)^2.
$$
\end{theorem}

\begin{theorem}\label{corol:main2}
The number of disjoint admissible domains in $\R{2}$ satisfies
$$
s(2)\leq 3.
$$
\end{theorem}

The proof of Theorem~\ref{th:main3} is given in
Section~\ref{sec:Li} below. Our argument uses the same idea as in
\cite{Li0} but we rearrange the original method in a more optimal
way to relax an exponential (in $n$) growth in (\ref{double}) to a polynomial growth.

The proof of Theorem~\ref{corol:main2} is more delicate, and it is
based on special bilateral estimates for the so-called angular density. More precisely, given a differentiable function $w(z)$ in an unbounded domain $G$, we define
the \textit{ angular density} of $G$ (with respect to $w$) as
$$
\Theta_w({G}):=\liminf_{R\to \infty}\frac{1}{\ln R}
\int\limits_{{G}_w(1,R)}\frac{\sqrt{1+|\nabla w|^2}}{(|x|^2+w^2(x)
)^{n/2}}dx,
$$
where
$
{G}_w(r,R)=\{x\in G: r^2<|x|^2+w^2(x)<R^2\}.
$
For $w\equiv 0$ we obtain the angular density (sometimes called the logarithmic volume)
$$
\Theta_0(G):=\liminf_{R\to \infty}\frac{1}{\ln R}
\int\limits_{G_0(1,R)}\frac{dx}{|x|^n}
$$
which is a metric invariant of $G$. For a solid cone over a spherical domain $\Omega\subset \Sirc$, $\Theta_0 (G)$
coincides with the surface measure of $\Omega$ and, in general, $\Theta_0(G)\leq \omega_{n-1}$, where $\omega_{n-1}$ is the $(n-1)$-dimensional Lebesgue measure of the unit sphere in $ \Sirc\subset \R{n}$.  In the two-dimensional case $\Theta_0(G)$ is well-known in function theory as the logarithmic area of $G$ (introduced by Teichm\"uller in 1930s, see \cite[Ch. VI]{Wittich}). In general, $\Theta_w({G})$ is very related the modulus of the family of curves (the extremal length) joining a compact on $\mathcal{G}$ with infinity, and it  can be thought of as a generalized logarithmic volume of the graph (\ref{xn}) at infinity. We would like also to mention that this generalized logarithmic volume, but for immersed higher-dimensional  minimal surfaces, is an effective tool for estimating of the number of ends of the surfaces, see, e.g. \cite{T94}, \cite{Asue}.


In this setup, Theorem~\ref{corol:main2} follows from the following
two results which are interesting  in their own right.

\begin{theorem}\label{th:main2}
Let $\mathcal{G}$ be an $n$-dimensional minimal graph (\ref{xn}) in $\R{n+1}$
supported on a domain $G\subset \R{n}$. Then
\begin{equation}\label{eee}
\Theta_w({G})\leq h_n \cdot \Theta_0(G),
\end{equation}
where
\begin{equation}\label{hn}
h_{n}:=(n-1)\int\limits_0^{+\infty}\frac{d\tau}{\left(1+\tau^2\right)^{n/2}}
=\frac{\sqrt{\pi}\Gamma\left(\frac{n+1}{2}\right)}{\Gamma(\frac{n}{2})}.
\end{equation}
\end{theorem}

\begin{theorem}\label{th:main1}
Let $\mathcal{G}$ be a two-dimensional minimal graph in $\R{3}$
supported on a non-trivial domain $G$. Then
\begin{equation}\label{e-star}
\Theta_w({G})\geq \pi\end{equation}
with equality when
$\mathcal{G}$ is a half-plane.
\end{theorem}


Throughout this paper we use the following  notation
$$
\begin{array}{ll}
\ox &=(x_1,\ldots,x_{n})\in\R{n};\\
\xv &=(\ox ,\xn )\in\R{n+1}; \\
\Pi\equiv\R{n}\hspace{0.5cm} &=\{\xv\in\R{n+1}:\xn =0\} ;\\
\Ball (R)&=\{\xv\in\R{n+1}: |\xv|<R\}; \\
\Sph(R)&=\partial\Ball (R);\\
\Disk (R)&=\{\ox \in\R{n}:|\ox |<R\};\\
\Sirc(R)&=\partial \Disk (R),\\
X& : \R{n+1}\to\Pi \,\mbox{the orthogonal projection};\\
\end{array}
$$

\section{The angular density estimates}\label{sec:1}

\subsection{The main inequality}
Let $M$ be an $n$-dimensional Riemannian manifold and $\M=(M,\uv)$
be a minimal hypersurface given by a proper isometric immersion
$\uv:M\to \R{n+1}$. We make no distinction between a
point $y\in M$ and its image $u(y)\in\M$. If the boundary $\partial M$ is
non-empty, it will be always assumed that
\begin{equation}\label{equ:image}
\uv(\partial M)\subset \Pi=\{\xv\in\R{n+1}:\xn=0\}.
\end{equation}

 Consider the following auxiliary function
$$
f(\xv):= |x|^{1-n}\varphi\left(\frac{\xn }{|x|}\right),
$$
where $
\varphi(t):=\int\limits_0^{t}
(1+\tau^2)^{-n/2}\,d\tau.
$

\begin{lemma}\label{lem:1}
Let $\M$ be a properly immersed minimal hypersurface in $\R{n+1}$.  Then
\begin{equation}\label{e:uv}
\frac{1}{|\uv|^n}\leq {h_n}\frac{|\scal{e_{n+1}}{
N}|}{|X(\uv)|^n}+ \scal{\nabla f(\uv)} {\nabla u_{n+1}},
\end{equation}
where $h_n=\varphi(+\infty)$ is defined by (\ref{hn}) and
$N$ is the unit normal field to $\M$.
\end{lemma}

\begin{proof}
We use the standard formalism of covariant differentiation on
immersed manifolds, see e.g. \cite{KN}. Then we have for the
gradients
\begin{equation}\label{equ:prem}
\nabla u_{n+1}=e_{n+1}^\top, \qquad \nabla |X(\uv)|=
\frac{X(\uv)^\top}{|X(\uv)|},
\end{equation}
where the symbols $\top$ and $\bot$ denote the projections on the
tangent space and normal spaces to $\M$ respectively. Denote
$$
\xi:=\frac{u_{n+1}}{|X(\uv)|}.
$$
Then
\begin{eqnarray*}
\nabla f(\uv)&=&-\biggl(\frac{(n-1)\varphi(\xi)}{|X(\uv)|^n}
+\frac{u_{n+1}\varphi'(\xi)}{|X(\uv)|^{n+1}} \biggr)\nabla |X(\uv)|+
\frac{\varphi'(\xi)}{|X(\uv)|^n}\nabla u_{n+1} =\\
&=&-\biggl(\frac{(n-1)\varphi(\xi)}{|X(\uv)|^{n-1}}+
\frac{u_{n+1}}{\left(|X(\uv)|^2+u_{n+1}^2\right)^{n/2}}\biggr)\frac{\nabla |X(\uv)|}{|X(\uv)|}
+\frac{\nabla u_{n+1}}{\left(|X(\uv)|^2+ u_{n+1}^2\right)^{n/2}},
\end{eqnarray*}
which by virtue of (\ref{equ:prem}) yields
$$
\nabla f(\uv)=-\frac{X(\uv)^\top}{|X(\uv)|^{n+1}} H(\xi)
+\frac{e_{n+1}^\top}{\left(|X(\uv)|^2+u_{n+1}^2\right)^{n/2}},
$$
with
$
H(t)=(n-1)\varphi_{} \left(t  \right)
+\frac{t}{\left(1+t^2\right)^{n/2}}.
$
Hence we have
\begin{equation}
\scal{\nabla f_{}(\uv)}{\nabla
u_{n+1}}=-\frac{\scal{X(\uv)^\top}{e_{n+1}^\top}} {|X(\uv)|^{n+1}}H(\xi)
+\frac{|e_{n+1}^\top|^2}{\left(|X(\uv)|^2+u_{n+1}^2\right)^{n/2}}.
\label{tkaequ2}
\end{equation}
Since  $e_{n+1}$ and $X(\uv)$ are mutual orthogonal as the vector
fields  in $\R{n+1}$, we infer
$$
\scal{X(\uv)^\top}{e_{n+1}^\top}=\scal{e_{n+1}-e_{n+1}^{\bot}}
{X(\uv)-X(\uv)^{ \bot}}=  -\scal{e_{n+1}}{ N}\scal{X(\uv)}{ N}.
$$
Similarly,
$
|e_{n+1}^\top|^2=1-\scal{e_{n+1}}{ N}^2.
$
Inserting the obtained relations in (\ref{tkaequ2}) yields
\begin{eqnarray*}
\scal{\nabla f(\uv)}{\nabla u_{n+1}}&=&\frac{1-\scal{e_{n+1}}{
N}^2} {\left(|X(\uv)|^2+u_{n+1}^2\right)^{n/2}}+\frac{\scal{X(\uv)}{
N}\scal{e_{n+1}}{ N}}{|X(\uv)|^{n+1}} H(\xi)=
\\
&=&\frac{\scal{e_{n+1}}{ N}}{|X(\uv)|^n}\left(\frac{\scal{X(\uv)}{ N}}{|X(\uv)|}
H(\xi)-\frac{\scal{e_{n+1}}{ N}{}}{\left(1+\xi^2\right)^{n/2}}\right)\\
&+&\frac{1}{\left(|X(\uv)|^2+u_{n+1}^2\right)^{n/2}}.
\end{eqnarray*}
Hence  we get the following identity
\begin{equation}
\label{H}
\scal{\nabla f(\uv)}{\nabla
u_{n+1}}=\frac{1}{|\uv|^n}
+\frac{\scal{e_{n+1}}{ N}}{|X(\uv)|^n}
\left(\frac{\scal{X(\uv)}{N}}{|X(\uv)|}H(\xi)-\frac{\scal{e_{n+1}}{N}}
{\left(1+\xi^2\right)^{n/2}}\right).
\end{equation}
On the other hand, the mutual orthogonality of $X(\uv)$ and
$e_{n+1}$ yields
$$
\scal{\frac{1}{|X(\uv)|}X(\uv)}{N}^2+\scal{e_{n+1}}{N}^2\leq 1,
$$
and applying Cauchy's inequality to (\ref{H}) we arrive at
\begin{equation}\label{equ:min}
\frac{1}{|\uv|^n}\leq\scal{\nabla f(\uv)}{\nabla u_{n+1}}
+\frac{|\scal{e_{n+1}}{ N}|}{|X(\uv)|^n} \Phi(\xi).
\end{equation}
Here
$$
\Phi(\xi)=\sqrt{\frac{1}{\left(1+\xi^2\right)^n}+\left((n-1) {\varphi
\left(\xi\right)}+\frac{\xi}{\left(1+\xi^2\right)^{n/2}}\right)^2}.
$$
It remains only to estimate the latter expression. To this aim, we
observe that $\Phi(\xi)$ is an even function, hence
$
\sup_{\xi\in\R{}}\Phi(\xi)=\sup_{\xi\geq 0}\Phi(\xi).
$
On the other hand, for $\xi\geq 0$
\begin{eqnarray*}
\Phi(\xi)\Phi'(\xi)&=&
\frac{n}{\left(1+\xi^2\right)^{{n/2}+1}}\left[(n-1){\varphi\left(\xi\right)}+\xi{\varphi'\left(\xi\right)}\right]
-\frac{n\xi}{\left(1+\xi^2\right)^{n+1}}=\\
&=&
\frac{n(n-1){\varphi\left(\xi\right)}}{\left(1+\xi^2\right)^{n+2/2}}\geq
0.
\end{eqnarray*}
 Thus $\Phi(\xi)$ is  increasing in $(0,+\infty)$, and
 it follows that
 $$
\sup_{\xi\in\R{}}\Phi(\xi)=\lim_{\xi\to\infty}
\Phi(\xi)=(n-1)\varphi(+\infty)=h_n.
 $$
Combining this with (\ref{equ:min}) proves the lemma.
\end{proof}

\subsection{Minimal quasigraphs}
A minimal hypersurface $\M$  is said to be a \textit{quasigraph}
if the orthogonal projection
$$
X\circ \uv: M\to \Pi\equiv \R{n}
$$
is a proper mapping. By $q(x)$, $x\in\Pi$, we denote the
multiplicity of the projection, i.e. $q(y)$ equals the number of
points $y\in M$ such that $X(\uv(y))=x$. Define the average
multiplicity as
$$
Q(t):=\frac{1}{ t^{n-1}}\int\limits_{\Sph(t)}q(x),
$$
where the integration is taken over the standard spherical
measure. Here and in  what follows, we suppress the notation of
differentials when integrals are taken over submanifolds. We have the following estimate (for $n=2$, a similar result was proven by  V.~Miklyukov and the author in \cite{MT87}).

\begin{lemma}\label{lem:s}
Let $\M$ be a minimal quasigraph. Then for sufficient large $r>0$,
and for all $R>r$
\begin{equation}\label{new-1}
\int\limits_{M(r,R)}\frac{1}{|\uv|^n}\leq c+
\frac{h_nQ(R)}{n-1}+h_n\int\limits_{r}^R\frac{Q(t)dt}{t},\end{equation}
where $M(r,R)=\{y\in M: r<|\uv(y)|<R\}$, and $c$ is some constant
depending on $r$.
\end{lemma}

\begin{proof}
We make use the notation of Lemma~\ref{lem:1}. If $\partial M\ne
\emptyset$ we put
$$
r_0:=\min _{y\in\partial M} |X( u(y))|,
$$
and $r_0=1$ otherwise. Define for $r_0<r<R$
$$
M^*(r,R)=\{y\in M: r<|X( u(y))|<R\}.
$$
By
 (\ref{equ:image}) we have
$ \varphi(\xi(y))=0$ for $y\in\partial M$, where $\xi=u_{n+1}/|X(\uv)|$.
Since $\M$ is minimal, the coordinate function $u_{n+1}$ is harmonic, hence we obtain by
integrating by parts
\begin{eqnarray*}
\int\limits_{M^*(r,R)}\scal{\nabla f(\uv)}{\nabla
u_{n+1}}&=&\int\limits_{\partial M^*(r,R)}f(\uv)\scal{\nabla
u_{n+1}}{\nu} \\
&=&\int\limits_{\mathcal{C}(r)\cup
\mathcal{C}(R)}\frac{\varphi(\xi)}{|X(\uv)|^{n-1}}\scal{e_{n+1}^\top}{\nu},
\end{eqnarray*}
where $ \mathcal{C}(t):=\{y\in M:|X(u(y))|=t\}, $ and $ \nu$ denotes
the unit normal field to $\mathcal{C}(t)$. Since
$|\varphi(\xi)|\leq h_{n}/(n-1)$, we obtain
\begin{equation}\label{equ:lem}
\biggl|\int\limits_{M^*(r,R)}\scal{\nabla f(\uv)}{\nabla u_{n+1}}
\biggr| \leq \frac{h_n}{n-1}( I(r)+I(R)),
\end{equation}
where
$$
I(t)=\frac{1}{t^{n-1}}\int\limits_{\mathcal{C}(t)}|\scal{e_{n+1}^\top}{\nu}|
$$

In order to estimate $I(t)$ we fix arbitrarily a regular value $t$
of the function $|X(u)|$ and consider the orthogonal
projection
$$
X^*:=X\circ\uv|_{\mathcal{C}(t)} : \mathcal{C}(t)\to \Sph(t).
$$
To derive the Jacobian of this mapping, we identify in a standard
way the tangent space $T_y\mathcal{C}(t)$ with its image in
$\R{n+1}$. It is not difficult to see that the Jacobian is found by the formula
\begin{equation}\label{Jac}
|\det d_yX^*|=|e_{n+1}^L|,
\end{equation}
where $L=L_y$ is the two-dimensional orthogonal complement to
$T_y\mathcal{C}(t)$ in $\R{n+1}$. Indeed,
choose an orthonormal basis $\{E_j\}_{j=1}^{n-1}$ of the tangent
space $T_y\mathcal{C}(t)$ such that $E_j\in\Pi$ for $1\leq j\leq
n-2$. Then
$$
d_yX^*(E_j)=d(X\circ\uv
)(E_j)=X({E}_j)={E}_j-e_{n+1}\scal{e_{n+1}}{{E}_j}.
$$
In particular, $d_yX^*(E_j)=E_j$ for $1\leq j\leq n-2$. Hence the (absolute value of)
Jacobian is found by
$$
|\det d_yX^*|=|X(E_{n-1})|=\sqrt{1-\scal{e_{n+1}}{{E}_{n-1}}^2}
$$
Since $\scal{e_{n+1}}{E_j}=0$ for all $1\leq j\leq n-2$, we get
(\ref{Jac}).

Next, using the orthogonality of $\nu$ and $T_y\mathcal{C}(t)$, we
conclude $\nu\in L_y$, and therefore (\ref{Jac}) yields
$$
|\det d_yX^*|\geq |\scal{e_{n+1}}{\nu}|.
$$
Applying the last inequality and the change variables formula, we
obtain
$$
I(t)\leq \frac{1}{t^{n-1}}\int\limits_{\Sph(t)}q(x)=Q(t),
$$
which,
in view of (\ref{equ:lem}), implies
\begin{equation}\label{sq}
\biggl|\int\limits_{M^*(r,R)}\scal{\nabla f(\uv)}{\nabla u_{n+1}}
\biggr| \leq \frac{h_n}{n-1}( Q(r)+Q(R)).
\end{equation}

Now notice that
$$
M(r,R)\subset M^*(r,R)\setminus {M(0,r)}\subset M^*(r,R)\cup K
$$
where $K=\overline{M^*(0,r)}\setminus M(0,r)$ is compact. Hence, we find from
(\ref{e:uv}) and (\ref{sq})
\begin{equation}
\begin{split}
\label{s1}
 \int\limits_{M(r,R)}\frac{1}{|\uv|^n}&\leq
\int\limits_{K}\frac{1}{|\uv|^n}+\int\limits_{M^*(r,R)}\frac{1}{|\uv|^n}
\leq \int\limits_{K}\frac{1}{|\uv|^n}
+\int\limits_{M^*(r,R)}
h_n\frac{|\scal{e_{n+1}}{
N}|}{|X(u)|^n}\\
&+\int\limits_{M^*(r,R)}
\scal{\nabla f} {\nabla u_{n+1}}\leq c+\frac{h_nQ(R)}{n-1}+{h_n}\int\limits_{M^*(r,R)}
\frac{|\scal{e_{n+1}}{ N}|}{|X(u)|^n} ,
\end{split}
\end{equation}
where
$
c=\frac{h_n}{n-1}Q(r)+\int\limits_{K}\frac{ 1}{|\uv|^n}
$
does not depend on $R$.

Similarly  to the proof of (\ref{Jac}), one can show that the
Jacobian of projection $X\circ\uv : M \to \Pi$ is equal to
$\scal{e_{n+1}}{N}$. Therefore, we have
\begin{equation}\label{equ:est}
\int\limits_{M^*(r,R)}\frac{|\scal{e_{n+1}}{N}|}{|X(u)|^n} = \int\limits_
{\Disk (R)\setminus \Disk (r)}\frac{q(x)}{|\ox |^{n}}d\ox
=\int\limits_{r}^R\frac{Q(t)}{t}dt.
\end{equation}
Thus, combining of (\ref{s1}) and (\ref{equ:est}) we get the
desired inequality.
\end{proof}

\subsection{Proof of Theorem~\ref{th:main2}} Now, let $\mathcal{G}$
be the minimal surface given as a graph (\ref{xn}) supported on a domain $G\subset \Pi$. Here $M=G$ and $ \uv(x)=(x,w(x)). $ Then the counting
function of the graph $\mathcal{G}$ coincides with the
characteristic function of $G$: $q(\ox)=\chi_{G}(\ox)$.
Hence,
$$
Q(t)=\frac{1}{ t^{n-1}}\int\limits_{\Sirc(t)}\chi_{G}(\ox)\leq \mathrm{Area}(\Sirc(1))=
\omega_{n-1}.
$$
Since $ G(r,R)$ is contained (up to a compact set) in $G\cap \Disk (r,R)$ for any $r>r_0$, where  $r_0$ is defined as in Lemma~\ref{lem:s}, we obtain from (\ref{new-1}) that
$$\int\limits_{G(r,R)}\frac{\sqrt{1+|\nabla w|^2}}{(|x|^2+w^2(x)
)^{n/2}}dx= \int\limits_{M(r,R)}\frac{1}{|\uv|^n}
\leq c+\frac{h_n\omega_{n-1}}{n-1}+h_n \int\limits_{G\cap (\Disk (R)\setminus \Disk (r))} \frac{d\ox}{ |\ox|^{n}},
$$
which easily implies (\ref{eee}), and the theorem is proved.

\section{The Dirichlet integral estimates}

\subsection{The weighted fundamental frequency}

In order to get lower estimates for $\Theta_w(G)$ we use a variation of the weighted fundamental frequency technique developed by V.~Miklyukov and the author in \cite{MT94}.
The origin of the method  gives rise to the classical Ahlfors distortion theorem on the Denjoy conjecture
relating the number of asymptotic values for an entire function and its order \cite{Baern}. In the early 1980s the method was used  by Miklyukov (see, e.g., \cite{M80}) in connection with upper estimates on the number of sublevel sets of solutions to a wide class of quasilinear PDE's in  higher dimensional Euclidean spaces.  Below we describe briefly some necessary definitions and facts (see also \cite{MT94}).

Let $\Sigma$ be a finite collection of (connected) one-dimensional
compact Riemannian manifolds with non-empty boundaries, and $g(y)$
be a smooth positive function defined on $\Sigma$.
Consider the following variational problem
\begin{equation}\label{e:lambda}
\lambda(\Sigma,g)=\inf_{\varphi}\left(\frac{\int\limits_{\Sigma}|D\varphi(y)|^2{g(y)^{-1}}\dsi}
{\int\limits_{\Sigma}\varphi^2(y)g(y)\dsi}\right)^{1/2}
\end{equation}
where $D\varphi$ stands for the covariant derivative of $\varphi$
with respect to the inner metric on $\Sigma$, and the infimum is
taken over all Lipschitz functions $\varphi(y)$ subject to the condition
$\varphi(y)=0$ on $\partial \Sigma$. The quantity $\lambda(\Sigma,g)$ is
called the \textit{weighted fundamental frequency} of $\Sigma$
(with respect to the weight $g$).

\begin{lemma}\label{lem:one-dim}
In the above notation,
$$
\lambda(\Sigma,g)\geq
\pi\biggl(\int\limits_{\Sigma}g(y)\dsi\biggr)^{-1}.
$$
\end{lemma}

\begin{proof}
We first assume that $\Sigma$ consists of a single component. Then
$\Sigma$ is isometric to certain Euclidean interval $I=[0,\beta]$,
where $\beta$ is the length of $\Sigma$. Denote by
$
f(t):I\to \Sigma
$
the corresponding isometry. Let $\varphi(y)$ be an arbitrary
Lipschitz function on $\Sigma$ subject to the zero Dirichlet boundary condition $\varphi|_{\partial \Sigma}=0$. Define $\psi(t)=\varphi\circ f(t)$, $G(t)=g\circ f(t)$. In this
notation we have
\begin{equation}\label{e:1}
\frac{\displaystyle\int_{\Sigma}|D\varphi(y)|^2{g(y)^{-1}}\dsi}
{\displaystyle\int_{\Sigma}\varphi^2(y)g(y)\dsi}=
\frac{\displaystyle\int_{0}^\beta \psi'^2(t)G(t)^{-1}dt}
{\displaystyle\int_{0}^\beta \psi^2(t)G(t)dt}.
\end{equation}
Let
$
\tau(t)=\int\limits_{0}^t G(\xi)d\xi$.
Since $g>0$, $\tau(t)$ is an increasing function. A new function $\zeta$ defined by $\psi(t)=\zeta(\tau(t))$ is obviously Lipschitz on $[0,\tau(\beta)]$ and satisfies $\zeta(0)=\zeta(\tau(\beta))=0$. It
follows from (\ref{e:1}) and Wirtinger's inequality \cite{HLP}
that
$$
\lambda^2(\Sigma,g)=\inf_{\zeta(0)=\zeta(\tau(\beta))=0}\frac{\displaystyle\int_{0}^{\tau(\beta)}
\zeta'^2(\tau)d\tau} {\displaystyle\int_{0}^{\tau(\beta)}
\zeta^2(\tau)d\tau}={\left(\frac{\pi}{ \tau(\beta)}\right)}^{2}
=\frac{\pi^2}{\displaystyle\left(\int_{\Sigma}g(y)\dsi\right)^2}.
$$

Returning to the general case, let $\Sigma=\cup_{j=1}^p \Sigma_j$ be decomposition of $\Sigma$ into a finite union of
connected components. It suffices only to show that
\begin{equation}\label{e:connect}
\lambda(\Sigma,g)=\min_{1\leq j\leq p} \lambda(\Sigma_j,g).
\end{equation}
Notice that the upper  bound
$
\lambda(\Sigma,g)\leq \min_{1\leq j\leq p} \lambda(\Sigma_j,g)
$
follows easily from the definition. On the other hand, let $k$ be an index for which
$$
\lambda(\Sigma_k,g)=\min_{1\leq j\leq p} \lambda(\Sigma_j,g)
$$
and let a Lipschitz function $\varphi_j$, $1\leq j\leq p$, be chosen arbitrarily such that
$\varphi_j(y)=0$ on $\partial \Sigma_j$.
Denote by $\varphi(y)$ the function on $\Sigma$ such that
$
\varphi(y)=\varphi_j(y)$ for $y\in\Sigma_j.$
Then for all $1\leq j\leq p$
$$
\int\limits_{\Sigma_j}|D\varphi_j(y)|^2{g(y)^{-1}}\dsi \geq
\lambda^2(\Sigma_k,g)
\int\limits_{\Sigma_j}\varphi_j(y)^2{g(y)}\dsi,
$$
which yields
$$
\int\limits_{\Sigma}|D\varphi(y)|^2{g(y)^{-1}}\dsi \geq
\lambda^2(\Sigma_k,g) \int\limits_{\Sigma}\varphi(y)^2{g(y)}\dsi.
$$
It follows from (\ref{e:lambda}) that $\lambda(\Sigma,g)\geq
\lambda(\Sigma_k,g)$, and (\ref{e:connect}) is proved.

\end{proof}

\subsection{The Dirichlet integral estimates}

Let us now describe our strategy of proof of Theorem~\ref{th:main1} in some more detail. A key auxiliary result is a lower energy estimate (\ref{e:main}) below. In order to establish it, we use a standard technique (sometimes called the Saint-Venant principle) of differentiating  the Dirichlet integrals taken over sublevel sets of a certain exhausting function (in our case, the distant function) and subsequent transforming of the obtained  integrals by using the weighted fundamental frequency (\ref{e:lambda}) into a differential inequality. Another important property, which ensures a.e.  differentiability of the Dirichlet integrals, is  monotonic character of the family of the sublevel sets.  In our case, however, there is an obstacle for a straightforward  using of the fundamental frequency, namely, level sets of the distant function can be a priori contain closed components. In order to get round the difficulty we cut off the `bad' components of sub-level sets in such a way that the remaining set still join the monotonic property. This preparation work is given in Lemmas~\ref{lem-G}, \ref{monosrt} below.

Let $\mathcal{G}=(G,w)$ be a  graph of a solution $w(x_1,x_2)$  of (\ref{mineq}) with the boundary condition (\ref{e:supp}), where  $G\subset \Pi$ is an admissible non-trivial domain and
$$
\Pi:=\{\overline{x}=(x_1,x_2,x_3)\in\R{3}:x_3=0\}.
$$

Denote by $G^-$ the union of all unbounded components of $\Pi\setminus \overline{G}$ (it is non-empty by non-triviality assumption about $G$), and set $G^+=\Pi\setminus \overline{G^-}$. Clearly $G\subset G^+.$

Denote by
$
\rho(x)=\sqrt{x_1^2+x_2^2+w^2(x_1,x_2)}
$
the distant function on $\mathcal{G}$ and consider the set
$$
G(t)=\{x\in G: \rho(x)<t\}.
$$
It follows immediately from the definition that
\begin{equation}\label{monott}
G(t_1)\subseteq G(t_2), \qquad t_1< t_2.
\end{equation}

We shall without loss of generality assume that $w(x_1,x_2)$ is a non-trivial solution, because for $w\equiv 0$ one has $\Theta_w({G})=\Theta_0({G})=\pi$. Since $w$ is non-trivial, it is not difficult to see (by virtue of (\ref{e:supp}) and the maximum principle) that for any  regular value  $t>0$ of the distance function $\rho(x)$ the \textit{relative boundary} $\partial G(t)\setminus \partial G$ is non-empty. Hence it splits into a finite collection of one-dimensional regular curves. A closed component of
$\partial G(t)\setminus \partial G$  will be called a \textit{cycle}. Obviously, any cycle $\Gamma$ is contained in the interior of $G$, in particular, $w(x)>0$ holds everywhere on $\Gamma$.

The remaining components of $\partial G(t)\setminus \partial G$ will be called \textit{arcs}. Any arc is contained in $G$ with the end-points on $\partial G$. It follows that function $w(x)$ vanishes on the boundary of any arc. An arc with both end-points in $G^+$ is called an \textit{inner} arc, otherwise it is called \textit{exterior}.

\begin{remark}\label{remm}
Any cycle or arc can be a part of the boundary of only one component of $G(t)$. Indeed, $\rho(x)\equiv t$ on $\partial G(t)\setminus \partial G$. It is well-known that $\rho(x)$ is subharmonic (in the metric of a minimal surface), hence it satisfies the strong maximum principle. It follows that $\rho(x)-t$ changes its sign in any neighborhood of any point from $\partial G(t)\setminus \partial G$, and the claim follows. For any curve $\gamma\in \partial G(t)\setminus \partial G$ we shall denote by $\mathcal{O}(\gamma)$ the unique component of $G(t)$ whose boundary contains $\gamma$.
\end{remark}

\begin{lemma}\label{lem-G}
Let $\Gamma$ be a cycle corresponding to some regular value $t>0$. Let $\mathcal{O}_\Gamma$ be the corresponding component of $G(t)$. Then

\begin{itemize}
  \item[(i)] $\mathcal{O}_\Gamma$ lies inside of $\Gamma$, i.e. $\mathcal{O}_G$ is a subset of the bounded component of $\Pi\setminus \Gamma$;
  \item[(ii)] $\overline{\mathcal{O}_\Gamma}\subset G^+$;
  \item[(i)] $\partial\mathcal{O}_\Gamma\setminus\Gamma\subset\partial G$.
  \end{itemize}
In particular, for any  component of $G(t)$, its  relative boundary consists of either  cycles or arcs.
\end{lemma}

\begin{proof}
Since $\Gamma$  is a closed curve, the complement $\Pi\setminus \Gamma$ splits into two components by the Jordan curve theorem. Denote by $U$ the bounded component of the complement. Then $\Gamma=\partial U$. Since $\Gamma$ is contained in $G$ with some neighborhood, the common part of $U$ and $G$ is non-empty and connected:
$$
U_G:=U\cap G\ne \emptyset.
$$
We claim that
\begin{equation}\label{max}
\rho(x)<t, \qquad \forall
x\in{U_G}.
\end{equation}
Since the distant function  $\rho(x)$ is subharmonic and non-constant, by the strong maximum principle
$$
\rho(x)<\max_{x\in\partial U_G}\rho(x) \qquad \forall
x\in{U_G},
$$
hence it suffices to show that $\rho(x)\leq t$ on $\partial U_G$. To this aim, we note that $\partial U_G$ consists of the `outer' part $\Gamma$ and the `inner' part which, if non-empty, is a subset of $\partial G$. On $\Gamma$  we   have trivially $\rho(x)= t$. Consider some point on the remaining part: $z_0\in \partial G\cap U$. Then $w(z_0)=0$ and $z_0$ is an interior point of
$U$. Since $|x|$ is subharmonic in $\Pi$, we get
$$
\rho(z_0)=|z_0|< \max _{x\in \overline{U}}|x|=\max _{x\in
\Gamma}|x|\leq\max _{x\in \Gamma}|\rho(x)|=t,
$$
hence $\rho(z_0)<t$. This proves (\ref{max}).

Applying (\ref{max}) and the strong maximum principle to $\rho(x)$ along $\Gamma$, it is easy to see that $\mathcal{O}_\Gamma\equiv U_G$ (see also Remark~\ref{remm}), which proves (i). Then (ii) follows from the obvious observation that $\Gamma$ can enclose only bounded components of $\Pi\setminus \overline{G}$. The  statement  (iii) is an easy corollary of inclusion $\partial U_G\setminus \Gamma\subset \partial G$.

Finally, notice that we have also proved that $\partial \mathcal{O}_\Gamma$ does not contain any arc (in fact, the relative boundary  of $\partial \mathcal{O}_\Gamma$ is precisely $\Gamma$). This finishes the proof.

\end{proof}

\begin{corollary}\label{corr1}
If $G$ is a non-trivial domain then there is $t_0(G)>0$ such that for any $t\geq t_0(G)$, $\partial G(t)$ contains at least one exterior arc.
\end{corollary}

\begin{proof}
Notice that the union  $G^-$ of all unbounded components of $\Pi\setminus \overline{G}$ is non-empty and choose $z_0\in \partial G^-$ arbitrarily.
Then $z_0\in \partial G^-\subset \partial G$ and  for any $t>t_0(G):= \rho(z_0)$ we have $z_0\in \overline{G(t})$. Suppose now that for some $t>t_0(G)$, the  boundary  $\partial G(t)$ does not contain an exterior arc, hence $\partial G(t)\setminus \partial G$ consists of only   cycles and interior arcs. It follows then from Lemma~\ref{lem-G} and the definition of an interior arc that $\overline{G(t)}\subset G^+$ which contradicts the choice of $z_0$.

\end{proof}

Now we are ready to present the mentioned in the beginning of the section monotone family. Consider any regular (for $\rho\circ u$) value $t\geq t_0(G)$, where $t_0(G)$ is chosen as in Corollary~\ref{corr1} and denote by $G_\alpha(t)$ the union of the components of $G(t)$ which are contained with their closures in $G^+$ (equivalently, $G_\alpha(t)$ is the union of those components of $G(t)$ whose  relative boundaries consists of  cycles or inner arcs). Set
$$
G_\beta(t)=G(t)\setminus G_\alpha(t).
$$

Summarizing Lemma~\ref{corr1} and Corollary~\ref{corr1}, we have the following.

\begin{lemma}\label{monosrt}
$G_\beta(t)$ is non-empty for any $t\geq t_0(G)$ and its relative boundary does not contain cycles.
For any regular $t_2>t_1\geq t_0(G)$: $G_\beta(t_1)\subseteq G_\beta(t_2)$.
\end{lemma}

\begin{proof}
We briefly comment only the last assertion.
By virtue of (\ref{monott}), it suffices only to check that $G_\beta(t_1)\cap G_\alpha(t_2)=\emptyset$. This easily follows from the property of $\alpha$-components: $\overline{G_\alpha(t)}\subset G^+$, while for $\beta$-components we have $\overline{G_\beta(t)}\cap \partial G^-\ne \emptyset.$

%
 \end{proof}

Let $G$ be an arbitrary non-trivial domain.
For any regular (for $\rho$) $t\geq t_0(G)$  consider the  Dirichlet integral
\begin{equation}\label{frrom}
J(t)=\int\limits_{{\mathcal{G}}_\beta(t)}|\nabla u_3|^2,
\end{equation}
where $\nabla$ denotes the covariant derivative with respect to
the inner metric of $\mathcal{G}$. Here and in what follows we denote for brevity  $\mathcal{G}_\beta(t)=X^{-1}({G}_\beta(t))$ etc. , where $X:\mathcal{G}\to G$ is as usual the orthogonal projection. Observe that $J(t)$ is increasing for all regular $t$ and redefine it in a standard way by setting $J(t)=\sup \{J(s):s<t \text{ and $t$ is regular}\}$. Thus obtained function is lower semi-continuous and non-decreasing in $(t_0(G),+\infty)$.

Let the weight function in (\ref{e:lambda}) is chosen as follows:
\begin{equation}\label{equ:gradient}
g(y):=|\nabla \rho(y)|=\frac{|\uv^{\top}(y)|}{|\uv(y)|}\leq 1,
\end{equation}
where $u=(x_1,x_2,w(x_1,x_2))$. Here we hold the above notation $\rho$ for the restriction of the distant function $\sqrt{x_1^2+x_3^2+x_3^2}$ on $\mathcal{G}$. Notice also that $g(y)>0$ because $t$ ia a regular value.

\begin{lemma}\label{lem:Dirichlet}
Let $\mathcal{G}$ be an two-dimensional minimal graph supported on a
non-trivial domain $G$. Then for any  $t_2>t_1>t_0({G})$
\begin{equation}\label{e:main}
J(t_2)\geq J(t_1)\exp \left(2
\int\limits_{t_1}^{t_2}\lambda(\Sigma(t),g)dt\right),
\end{equation}
where $\Sigma(t)=\partial \mathcal{G}_\beta(t)\setminus \partial \mathcal{G}$ is the relative boundary of $\mathcal{G}_\beta(t)$ in $\mathcal{G}$
and $t_0(G)$ is chosen as in Corollary~\ref{corr1}.
\end{lemma}

\begin{proof}
By Lemma~\ref{monosrt}, $J(t)$ is an non-decreasing function, hence it is differentiable  almost everywhere. Hence, the set $T$ of all \textit{regular} points $t>t_0({G})$ where $J(t)$ is differentiable has full measure in $(t_0(G),+\infty)$. Let $t\in T$ be chosen arbitrarily.
By our construction, the  relative boundary  $\Sigma(t)$ splits in a finite collection of arcs with end-points on $\partial \mathcal{G}=\partial G$, hence
\begin{equation}\label{2}
u_3(y)=0, \quad y\in\partial \mathcal{G}.
\end{equation}
We have by (\ref{e:lambda})
\begin{equation}\label{equ:lam}
\int\limits_{\Sigma(t)}u_3^2g\dsi\leq
\frac{1}{\lambda^2(t)}
\int\limits_{\Sigma(t)}|Du_3|^2g^{-1}\dsi,
\end{equation}
where $\lambda(t)=\lambda(\Sigma(t),g)$ and $D$ denotes the induced on $\Sigma(t)$ covariant derivative (in this case, the directional derivative along the unit tangent vector to $\Sigma(t)$). On the other hand, by harmonicity of $u_3$ and (\ref{2}) we obtain
\begin{equation}\label{29}
J(t)=\int\limits_{\mathcal{G}_\beta(t)}|\nabla
u_3|^2=\int\limits_{\mathcal{G}_\beta(t)}\diver (u_3 \nabla u_3)=
\int\limits_{\partial \mathcal{G}_\beta(t)}u_3\scal{\nabla
u_3}{\nu}\dsi=\int\limits_{\Sigma(t)}u_3\scal{\nabla u_3}{\nu}\dsi,
\end{equation}
where $\nu$ stands for the unit outward  normal to $\Sigma(t)$. By the Cauchy inequality,
$$
|u_3\scal{\nabla u_3}{\nu}|\leq
\frac{\lambda(t)}{2}u_3^2g+\frac{1}{2\lambda(t)g} |\scal{\nabla
u_3}{\nu}|^2.
$$
and applying  (\ref{equ:lam}) and (\ref{29}) we obtain
\begin{equation}\label{equ:Jprime2}
\begin{split}
J(t)&=\int\limits_{\Sigma(t)}u_3\scal{\nabla u_3}{\nu}\dsi
\leq
\frac{\lambda(t)}{2}\int\limits_{\Sigma(t)}u_3^2g\dsi+\frac{1}{2\lambda(t)}\int\limits_{\Sigma(t)}
|\scal{\nabla u_3}{\nu}|^2g^{-1}\dsi\\
&\leq
\frac{1}{2\lambda(t)}\int\limits_{\Sigma(t)}(|Du_3|^2+|\scal{\nabla
u_3}{\nu}|^2)g^{-1}\dsi.
\end{split}
\end{equation}
Note that $|\nabla u_3|^2=|Du_3|^2+\scal{\nabla u_3}{\nu}^2.$
By our choice of $t$, $J(t)$ is differentiable at $t$, hence applying the co-area formula we find from (\ref{equ:Jprime2})
\begin{equation}\label{aff}
J(t)\leq \frac{1}{2\lambda(t)}\int\limits_{\Sigma(t)}
\frac{|\nabla
u_3|^2}{g}\dsi=\frac{1}{2\lambda(t)}\int\limits_{\Sigma(t)}
\frac{|\nabla u_3|^2}{|\nabla
\rho|}\dsi=\frac{J'(t)}{2\lambda(t)}.
\end{equation}

Thus the  differential inequality $\frac{d}{dt}\ln J(t)\geq 2\lambda(t)$ holds for almost all $t>t_0(G)$. In addition  $\ln J(t)$ is an non-decreasing function, hence for any $t_2>t_1>t_0({G})$:
$$
\ln J(t_2)-\ln J(t_1)\geq \int_{t_1}^{t_2}(\ln J(t))'\, dt=\int_{t_1}^{t_2}\frac{J'(t)}{J(t)}\,dt
$$
and after applying of (\ref{aff}), we get
$$
\ln J(t_2)-\ln J(t_1)\geq
2\int_{t_1}^{t_2}\lambda(t) \,dt
$$
and   (\ref{e:main}) follows.

\end{proof}


\subsection{Proof of Theorem~\ref{th:main1}}

We use notation of Lemma~\ref{lem:Dirichlet}. By
Cauchy's inequality we have for any $t_2>t_1>t_0(G)$
$$
\int\limits_{t_1}^{t_2}\lambda(\Sigma(t),g)dt \cdot
\int\limits_{t_1}^{t_2}\frac{1}{\lambda(\Sigma(t),g)}\frac{dt}{t^2}
\geq  (\ln t_2-\ln t_1)^2.
$$
Applying Lemma~\ref{lem:one-dim} we obtain
\begin{equation}\label{ff}
\int\limits_{t_1}^{t_2}\lambda(\Sigma(t),g)dt \geq
\frac{\displaystyle\pi (\ln t_2-\ln
t_1)^2}{\displaystyle\int\limits_{t_1}^{t_2}\frac{dt}{t^2}
\int\limits_{\Sigma(t)}g\dsi},
\end{equation}
where $g(y)$ is defined by (\ref{equ:gradient}). Thus, by
the co-area formula and inequality in (\ref{equ:gradient}), we get
$$
\int\limits_{t_1}^{t_2}\frac{dt}{t^2}
\int\limits_{\Sigma(t)}g\dsi=
\int\limits_{\mathcal{G}_\beta(t_1,t_2)}\frac{|\nabla
\rho|^2}{|\uv|^2} \leq
\int\limits_{\mathcal{G}_\beta(t_1,t_2)}\frac{1}{|\uv|^2},
$$
where $\mathcal{G}_\beta(t_1,t_2)=\mathcal{G}_\beta(t_2)\setminus\overline{\mathcal{G}_\beta(t_1)}$. Hence, combining the latter inequality with
(\ref{ff}) and (\ref{e:main}), we arrive at
\begin{equation*}\label{equ:two}
\ln \frac{J(t_2)}{J(t_1)} \geq2\int\limits_{t_1}^{t_2}\lambda(\Sigma(t),g)dt  \geq \frac{\displaystyle2\pi(\ln t_2-\ln
t_1)^2}{\displaystyle\int\limits_{\mathcal{G}_\beta(t_1,t_2)}\frac{1}{|\uv|^2}}.
\end{equation*}
or after rearrangement,
\begin{equation}\label{fol}
\frac{1}{\ln
\frac{t_2}{t_1}}\;\int\limits_{\mathcal{G}_\beta(t_1,t_2)}\frac{1}{|\uv|^2} \geq
\frac{2\pi(\ln t_2- \ln t_1)}{\ln J(t_2)- \ln J(t_1)}.
\end{equation}

On the other hand, $ |\nabla u_3|=|e_3^\top|\leq 1 $, hence we find from (\ref{frrom}): $
J(t)\leq \mathrm{Area}(\mathcal{G}_\beta(t))$. The
area growth estimate \cite[Lemma~1]{Li0} (see also
(\ref{equ:a-priori}) below) yields the quadratic area growth for minimal graphs:
$$ \mathrm{Area}(\mathcal{G}_\beta(t))\leq
 \mathrm{Area}(\mathcal{G}(t))\leq  3\pi R^2.
 $$
 Therefore,
$$
\frac{2\pi(\ln t_2- \ln t_1)}{\ln J(t_2)- \ln J(t_1)}\geq
\frac{2\pi(\ln t_2- \ln t_1)}{2\ln t_2 +\ln (3\pi)- \ln J(t_1)}
$$
 and it follows from (\ref{fol}) that
$$
\Theta_w({G})= \liminf_{t_2\to\infty}\frac{1}{ \ln
{t_2}}\;\int\limits_{\mathcal{G}(t_1,t_2)}\frac{1}{|\uv|^2} \geq
\liminf_{t_2\to\infty}\frac{1}{ \ln
{t_2}}\;\int\limits_{\mathcal{G}_\beta(t_1,t_2)}\frac{1}{|\uv|^2} \geq\pi,
$$
which proves (\ref{e-star}).

\subsection{Proof of Theorem~\ref{corol:main2}}\label{sec:2}

Let $\{G_j\}_{j=1}^{s}$ be an arbitrary admissible collection in
$\R{2}$. Without loss of generality we can assume that all $G_j$
are non-trivial, otherwise the maximum principle for solutions of
(\ref{mineq}) implies $s=1$. Then for $R>1$
$$
2\pi\ln R=\int\limits_{\Disk (1,R)}\frac{dx}{|x|^n}\geq \sum_{j=1}^s
\int\limits_{G_j\cap \Disk (1,R)}\frac{dx}{|x|^n},
$$
hence
$$
2\pi \geq \sum_{j=1}^s \Theta_0(G_j).
$$

On the other hand, as an immediate corollary of (\ref{eee}) and
Theorem~\ref{th:main2} for $n=2$ we have
$\Theta_0(G_j)\geq 2.$
Combining the obtained inequalities, we arrive at $ 2\pi \geq 2s $,
which finishes the proof.

\section{Proof of Theorem~\ref{th:main3}}\label{sec:Li}
Consider an arbitrary admissible collection of minimal graphs
$\mathcal{G}_j=(G_j,w_j)$, $j=1,\ldots,N$. Without loss of
generality, we may assume that $w_j(x)>0$ on $G_j$. As before, we
identify $\mathcal{G}_j$ with the corresponding submanifold in
$\R{n+1}$ equipped with the induced metric. Let
$\mathcal{G}_j(R)=\mathcal{G}_j\cap \Ball (R)$, and fix $r_0>0$ such
that all sets $\mathcal{G}_j(r_0)$ are non-empty. Denote by $W_j(y)$
the lifting of $w_j$ on the surface
$$
\mathcal{G}=\bigcup_{1\leq i\leq N}\mathcal{G}_i,
$$
i.e., $W_j(y)=w_j(\ox )$, if $y= (\ox ,w_j(\ox
))\in\mathcal{G}_j$, and $W_j(y)=0$ otherwise. By the disjointness
condition, $W_j(y)$ is a well-defined smooth function on
$\mathcal{G}$. Set
$$
\mathcal{G}(R):=\mathcal{G}\cap \Ball
(R)=\cup_{j=1}^{N}\mathcal{G}_j(R).
$$
and consider the following function
$$
W^R(y )=\sum_{j=1}^N\frac{W_j(y )}{\alpha_j(R)}, \qquad y
\in\mathcal{G},
$$
where $R>r_0$, and
\begin{equation}\label{e:max}
\alpha_j(R)=\max_{y\in\overline{\mathcal{G}_j(R)}}W_j(y)>0.
\end{equation}
It follows from the definition of $\alpha_j(R)$ that
$$
\max_{y\in\overline{\mathcal{G}(R)}}W^R(y )=\max_{1\leq j\leq
N}\;\max_{y\in\overline{\mathcal{G}(R)}}\frac{W_j(y
)}{\alpha_j(R)}=1
$$

Using the fact that all $\mathcal{G}_j$ are disjointly supported,
we obtain after integrating
\begin{equation}\label{equ:mean}
\sum_{j=1}^N
\frac{\Lambda(\mathcal{G}_j,R)}{\alpha_j(R)}=\int\limits_{\mathcal{G}(R)}W^R\leq
|\mathcal{G}(R)|,
\end{equation}
where $|\mathcal{G}(R)|$ denotes the volume of $\mathcal{G}(R)$, and
$$
\Lambda(\mathcal{G}_j,R):=\int\limits_{\mathcal{G}_j(R)}W_j.
$$

On the other hand, by Lemma~1 in \cite{Li0} we have
\begin{equation}\label{equ:a-priori}
|\mathcal{G}\cap \Ball (R)|\leq (n+1) |G\cap \Disk (R)|
\end{equation}
for all $R>0$. This implies
\begin{eqnarray}
|\mathcal{G}(R)|&=&\sum_{j=1}^{N}|\mathcal{G}_j(R)|\leq
(n+1)\sum_{j=1}^{N}|G_j\cap \Disk (R)|\leq \nonumber
\\ &\leq&  (n+1)
|\Disk (R)|=(n+1)\omega_nR^n,\nonumber
\end{eqnarray}
where $\Omega_n=|\Disk (1)|$ is the volume of the $n$-dimensional
unit ball.

For each $j$, $1\leq j\leq N$, there exists point
$y_j(R)\in\overline{\mathcal{G}_j(R)}$ at which the corresponding
maximum in (\ref{e:max}) is attained. Note that $W_j$ is a
harmonic function on $\mathcal{G}$. Hence, applying  the mean
value theorem \cite[Lemma~2]{Li0} for $W_j(y )$ and taking into
account the non-negativity of $W_j(y )$, we obtain for any $r>0$
$$
\alpha_j(R)=W_j(y_j(R))\leq
\frac{1}{\Omega_nr^n}\int\limits_{\mathcal{G}_j\cap \Ball (r;y_j(R))} W_j
\leq \frac{1}{\Omega_nr^n}\int\limits_{\mathcal{G}_j(r+R)}W_j=
\frac{\Lambda(\mathcal{G}_j,R+r)}{\Omega_nr^n},
$$
where $\Ball (r;y_j(R))=\{\overline{x}\in \R{n+1}:|\overline{x}-y_j(R)|<r\}$.
Then inserting the last inequalities in (\ref{equ:mean}) for $r>0$
and $R>r_0$ yields
$$
\sum_{j=1}^N
\frac{\Lambda(\mathcal{G}_j,R)}{\Lambda(\mathcal{G}_j,R+r)}\leq
(n+1)\frac{R^n}{r^n},
$$
or what is the same
\begin{equation}\label{last2}
\sum_{j=1}^N \frac{f_j(R)\;\;}{f_j(R+r)}\leq
(n+1)\left(\frac{R+r}{r}\right)^n,
\end{equation}
where $ f_j(t):=\Lambda(\mathcal{G}_j,t)t^{-n}. $

Let $\beta>1$ is given arbitrarily, and $r=(\beta-1)R$. Then using
the arithmetic-geometric means inequality in the left-hand side of
(\ref{last2}), we obtain
$$
(n+1)\left(\frac{\beta}{\beta-1}\right)^n\geq
N\left(\prod_{j=1}^{N}\frac{f_j(R)\;\;}{f_j(\beta
R)}\right)^{1/N}.
$$
Letting $R=\beta^k \rho$, $k=0,1,\ldots,m$, in the latter
inequality we get
\begin{equation}\label{last}
(n+1)\left(\frac{\beta}{\beta-1}\right)^n\geq
N\left(\prod_{j=1}^{N}\frac{f_j(\rho)\;\;}{f_j(\beta^m
\rho)}\right)^{1/Nm}
\end{equation}

On the other hand, $W_k(y )\leq t$ for  $y\in \mathcal{G}_k(t)$,
so we infer from (\ref{equ:a-priori})
\begin{equation}\label{sublinear}
f_j(t)=\frac{1}{t^n}\int\limits_{\mathcal{G}_j(t)}W_j\leq
\frac{|\mathcal{G}_j(t)|}{t^{n-1}}\leq (n+1)\Omega_n t.
\end{equation}
Hence (\ref{last}) and (\ref{sublinear}) yield
\begin{equation}\label{last1}
(n+1)\left(\frac{\beta}{\beta-1}\right)^n\geq N\frac{\lambda
^{1/Nm}}{\beta},
\end{equation}
where
$$
\lambda=\frac{1}{((n+1)\Omega_n\rho)^N}\prod_{j=1}^{N}f_j(\rho).
$$
Letting $m\to \infty$ in (\ref{last1}), we obtain
$$
N\leq (n+1)\beta\left(\frac{\beta}{\beta-1}\right)^n.
$$
The minimum of the latter right-hand side is attained at
$\beta=n+1$, and it follows that
$$
N\leq (n+1)^2\left(1+\frac{1}{n}\right)^n< e(n+1)^2
$$
and the theorem  follows.

\textbf{Acknowledgements}.
I am deeply indebted to the anonymous referee for many valuable comments on the previous version of the article.
I would like also to thank Yiu-Ming Lo and Jiaping Wang for several discussions  concerning the proof of theorem ~\ref{th:main3}.

\end{document}